\documentclass[12pt,a4paper]{article}

\usepackage[english]{babel}
\usepackage{amsmath}
\usepackage{amsfonts}
\usepackage{eufrak}
\usepackage{amssymb}
\usepackage{psfig}
\usepackage{epsf}
\usepackage{rotate}
\usepackage[all]{xy}
\usepackage{xypic}
\input xy
\xyoption{all}                 \begin{tiny}  \end{tiny}

  \newcommand{\St}{\mathfrak{S}}
 \newcommand{\K}{\mathrm{K}}
\newcommand{\Sp}{\mathrm{Sp}}

\newcommand{\la}{\lambda}

\newcommand{\C}{\mathrm{C}}
\newcommand{\M}{\mathcal{M}_{g,2}}

\newcommand{\A}{\mathcal{A}_{g,2}}

\newcommand{\PP}{\mathbb{P}}

 \newcommand{\J}{\mathcal{J}}
 \newcommand{\Z}{\mathbb{Z}}
\newtheorem{defi}{ Definition}[section]
\newtheorem{theo}[defi]{ Theorem}
\newtheorem{prop}[defi]{ Proposition}

\newtheorem{rema}[defi]{ Remark}
\newtheorem{lemm}[defi]{ Lemma}
 \newtheorem{ex}[defi]{ Example}

\begin{document}
 \setcounter{page}{1}

 \title{Vanishing thetanulls and hyperelliptic curves}
\author{Olivier Schneider}

 \maketitle

\begin{abstract}
Let $\mathcal{M}_{g}(2)$ be the moduli space of curves of genus $g$ with a level-2 structure.
 We prove  that there is always a non hyperelliptic element in the intersection of four thetanull divisors
  in $\mathcal{M}_{6}(2)$. We prove also that for all $g\geqslant3$, each component of the hyperelliptic locus in  $\mathcal{M}_{g}(2)$
  is a connected component of the intersection of $g-2$ thetanull divisors.
\end{abstract}

 \section{Introduction}

 \addcontentsline{toc}{section}{Introduction}
  Let $C$ be a curve of genus $g\geqslant3$ over $\mathbb{C}$ and $\sigma : \Z^{2g} \longrightarrow \mathrm{H}_{1}(\C,\Z)$
  a symplectic isomorphism. We associate to $(C,\sigma)$ a period matrix    $Z$ in the Siegel upper half space $\mathrm{H}_{g}$
  and the numbers
       $$ \theta[k](Z):=   \sum_{r \in \Z^{g}}
                    \exp \left(  \pi \sqrt{-1} \  \left[  ^{t}\left( r+\frac{k^{'}}{2} \right)
                               Z \left( r+\frac{k^{'}}{2} \right)  +  2 \  ^{t}\left( r+\frac{k^{'}}{2} \right)  \frac{k^{''}}{2}  \right] \right). $$
  for each  $k=(k^{'},k^{''})$ in  $(\Z/2\Z)^{2g}$ such as $k^{'}.k^{''}$ is even. The vanishing of $\theta[k](Z)$ depends
  only of $\sigma(\mathrm{mod} \  2)$, that is, on the class of $(C,\sigma)$ in $\mathcal{M}_{g}(2)$, the moduli space of curves of genus $g$ on
   $\mathbb{C}$ with a level-2 structure. Thus the zero locus of $ \theta[k]$ defines a divisor in
    $\mathcal{M}_{g}(2)$  called thetanull divisor of  characteristic $[k]$. In genus $3$, each thetanull divisor  is a component of the hyperelliptic locus
 in $\mathcal{M}_{3}(2)$. In genus $4$, we know since Riemann that each intersection of
 two thetanull divisors is
 an union of hyperelliptic components  in $\mathcal{M}_{4}(2)$.
  In genus  $5$,  Accola has established in  [\textbf{A1}], that with a condition on the three characteristics,
  the intersection  of the corresponding divisors
    is  an union of hyperelliptic components. We propose here to
 prove  that this fails in genus $6$ :
  \begin{theo}
   Each sub-variety of  $\mathcal{M}_{6}(2)$  intersection of four thetanull divisors
   contains an element wich is not hyperelliptic.      \label{f}
\end{theo}
  The first step of the proof is to classify the orbits of the action of  $\Sp_{2g}(\Z/2\Z)$ on the set of quadruplets
  of thetanull divisors. Afterwards,  we  finish the proof by verifying that one element of each orbit (and then every quadruplet) defines
  a subvariety of $\mathcal{M}_{6}(2)$ which contains a bi-elliptic curve. In view of this result, it is unlikely that hyperelliptic
   components are intersection of  $g-2$ thetanull divisors in higher genus.
 Nevertheless,  in the last part, we  prove the following result :
    \begin{theo}
   Every component of the hyperelliptic locus in $\mathcal{M}_{g}(2)$ can be defined as
   a connected component  of the intersection of $g-2$ thetanull divisors.  \label{trans}
   \end{theo}

   \section{Preliminaries}
  In what follows we denote by $\C$ a curve of genus $g$ over $\mathbb{C}$.
  Let  $\A$ be the moduli space of $g$-dimensional principally polarized abelian varieties with a level-2
  structure. It can be described as follows (see for instance [\textbf{BL}] chapter 4).
 Let   $\mathrm{H}_{g}$  be the Siegel generalized half-space :
  $$ \mathrm{H}_{g}=\{Z \in \mathrm{M}_{g}(\mathbb{C}) \ \mid \ \ ^{t}Z=Z \ ,\ \mathrm{Im} Z >0 \}.$$
  Then to each  $Z$ in $\mathrm{H}_{g}$ corresponds  the complex torus  $\mathbb{C}^{g}/(\Z^{g}+Z.\Z^{g})$ which
comes with a natural principal polarization and a level-2 structure. Moreover we
let $$\Gamma_{g}(2):=\{ M \in \Sp_{2g}(\Z) \ \mid \ M \equiv 1_{2g}   ( \mathrm{mod} \ 2 ) \}. $$ this group acts on
 $\mathrm{H}_{g}$ in the following way :
                           $$ M.Z:=(AZ+B)(CZ+D)^{-1} \ ,\ M=\left( \begin{array}{cc} A & B \\ C & D \end{array} \right) \in \Gamma_{g}(2).$$
    Then   $\A$ is isomorphic to the quotient space
  $\mathrm{H}_{g}/\Gamma_{g}(2)$. We denote by $\M$ the space
   which parametrizes the pairs $(\C,\sigma)$ where  $\C$ is  a curve of genus $g$ and
    $\sigma : (\Z/2\Z)^{2g} \longrightarrow \mathrm{H}_{1}(\C,\mathbb{Z}/2\mathbb{Z})$ a symplectic isomorphism.
    $\Sp_{2g}(\Z/2\Z)$ acts on $\M$ in a natural way. Let   $k=(k^{'},k^{''}) \in  (\Z/2\Z)^{2g}$ such as $k^{'}.k^{''}$ is even
    (in what follows, we will say that $k$ is even).
  Let
 $$\begin{array}{rccl}
            \theta[k]:  & \mathrm{H}_{g}  & \longrightarrow   &  \mathbb{C} \end{array}$$
  be the function defined by
  $$  \theta[k](Z)= \sum_{r \in \Z^{g}}
               \exp \left(  \pi \sqrt{-1} \  \left[  ^{t}\left( r+\frac{k^{'}}{2} \right)
                          Z \left( r+\frac{k^{'}}{2} \right)  +  2 \  ^{t}\left( r+\frac{k^{'}}{2} \right) \frac{k^{''}}{2}  \right] \right). $$

 Let $M$ be in $\Gamma_{g}(2)$. As a consequence of the transformation formula
 (see for instance [\textbf{I}] page 176), $\theta[k](M.Z)$ is proportionnal to $\theta[k](Z)$.
     By this fact, the zero locus of $\theta[k]$ defines a divisor of
     $\mathcal{M}_{g}(2)$ denoted by
        $\theta[k]_{0}$ and called the \textbf{thetanull divisor}  of characteristic $[k]$.
            Finally, the natural action of $\mathrm{Sp}_{2g}(\mathbb{Z}/2\mathbb{Z})$ on $(\mathbb{Z}/2\mathbb{Z})^{2g}$
            induces an action on
           these hypersurfaces : For each $k$ in $(\mathbb{Z}/2\mathbb{Z})^{2g}$ wich is even, for each $M$ in
           $\mathrm{Sp}_{2g}(\mathbb{Z}/2\mathbb{Z})$, we have

                   $$ M.\theta[k]_{0}=\theta[M.k]_{0} .$$

   \section{Theta characteristics, symplectic torsors}
                     Each concept of this part can be found (for instance) in [\textbf{S}].
         \begin{defi}
      Let  $\C$ be a curve of genus $g$ over $\mathbb{C}$ and $\K$ its canonical divisor class. A  \textbf{theta characteristic}
       $\Delta$ on  $\C$ is a degree $g-1$ divisor class such as  $$2\Delta \equiv \K$$
       Moreover  $\Delta$ is  called \textbf{even} (resp \textbf{odd}) if $h^{0}(\Delta) \equiv  0 \ ( \mathrm{mod} \ 2)$
        (resp $1 \ ( \mathrm{mod} \ 2 )$). \\
                \end{defi}

         Let   $\mathcal{S}(\C)$ be the set of theta characteristics on  $\C$ and let $\mathcal{S}^{+}(\C)$ be the set of even ones.
           \begin{defi}
       Let $(J,.)$  be a symplectic pair (that is a $\Z/2\Z$-vector space endowed with a non degenerate, alternate, bilinear form), we say that a pair $(S,Q)$
       is a symplectic torsor over $(J,.)$ if there is simply transitive action of $J$ on $S$ denoted $+$ and a
       mapping  $Q:S \longrightarrow \Z/2\Z$ having the property  \\
       $\forall s \in S$, $\forall j_{1} \in J$, $\forall j_{2} \in J$,
        $$Q(s)+Q(j_{1}+s)+Q(j_{2}+s)+Q(j_{1}+j_{2}+s)=j_{1}.j_{2}  \ \ \ \ (\ast)$$
    \end{defi}

     \begin{ex}
   \emph{Let $W$ be a set of order $2g+2$. We denote by  $P_{2}(W)$ the set of partitions on $W$ into two subsets.
     For each  subset $A$ of $W$, we denote by $\bar{A}$ the element $(A;\complement_{W}A)$ of $P_{2}(W)$.
     One has
     \begin{enumerate}
           \item An addition on $P_{2}(W)$ :
           For each $\bar{A}$ and $\bar{B}$ in $P_{2}(W)$
           $$\bar{A}+\bar{B}:=\overline{A \bigcup B-A \bigcap B}$$
           \item A map $p:P_{2}(W) \longrightarrow \Z/2\Z$ defined by
               $$p(\bar{A})=\mid A \mid ( \mathrm{mod} \ 2).$$
           We denote by $P_{2}^{+}(W)$ (resp $P_{2}^{-}(W)$) the set $p^{-1}(0)$ (resp $p^{-1}(1)$).
           \item  A map $e:P_{2}^{+}(W) \times P_{2}^{+}(W) \longrightarrow \Z/2\Z$ defined by
                    $$e(\bar{A};\bar{B}):=\mid A \bigcup B \mid  ( \mathrm{mod} \ 2).$$
          \end{enumerate}
     If $g$ is even (resp odd) then we have a quadratic form  defined on $P_{2}^{-}(W)$ (resp on  $P_{2}^{+}(W)$)  by
   $$Q_{-}(\bar{A})=\frac{\mid A \mid +1}{2} ( \mathrm{mod} \ 2)$$
   $$\left( resp  \ \ Q_{+}(\bar{A})=\frac{\mid A \mid }{2} ( \mathrm{mod} \ 2) \right). $$
   $(  P_{2}^{-}(W),Q_{-})$ (resp $(  P_{2}^{+}(W),Q_{+})$) is a symplectic torsor over $( P_{2}^{+}(W),e) $ denoted by $T(W)$.}
   \end{ex}

   \begin{ex}
   \emph{For each  symplectic pair  $(J,.)$, the set  of quadratic forms
     $$Quad(J):=\{q:J \longrightarrow \Z/2\Z \mid \forall j_{1} \in J \ , \ \forall j_{2} \in  J \ ,\
           q(j_{1})+q(j_{2})+q(j_{1}+j_{2})=j_{1}.j_{2}   \}$$
    endowed with the mapping  $\mathrm{Arf} :Quad(J)   \longrightarrow  \Z/2\Z$ (see for instance [\textbf{Sc}])
    is a symplectic torsor over $(J,.)$. Moreover we have the following property : $\forall q \in Quad(J) \ , \ \forall j \in J$,
 $$\mathrm{Arf}(j+q)=\mathrm{Arf}(q)+q(j)$$ }     \label{parity}
 \end{ex}

     \begin{prop}
      Let  $(\J_{2}(\C),.)$ be the  $\Z/2\Z$-vector space of order two points of the Jacobian of $\C$
      endowed with the intersection pairing. If we denote by
       $Q$  the mapping
      $$\begin{array}{rcl}
                    Q: \mathcal{S}(\C) & \longrightarrow & \Z/2\Z \\
                    \Delta    & \longmapsto       &  h^0(\Delta) \mod \ 2 \end{array}$$
     then $(\mathcal{S}(\C),Q)$ is a symplectic torsor over   $(\J_{2}(\C),.)$.
  \end{prop}

  \begin{rema}
  \emph{Let $\C$ be an hyperelliptic curve. If we denote by $W(\C)$ the set of its Weierstrass points, then we have a canonical
  isomorphism of symplectic torsors
  $$ \widetilde{\beta_{\C}}: T(W(\C)) \longrightarrow (\mathcal{S}(\C),Q)$$
   compatible with the canonical symplectic isomorphism $$ \beta_{\C}:(P_{2}^{+} (W(\C)),e) \longrightarrow (\J_{2}(\C),.).$$}
  \end{rema}
  Finally, we recall the result established by  Mumford in [\textbf{M}] :
  \begin{prop}
  The map
        $$\begin{array}{ccc}
               (\mathcal{S}(\C),Q) & \longrightarrow &   (Quad(\J_{2}(\C)), \mathrm{Arf})   \\
                \Delta  & \longmapsto & \left(  \begin{array}{rcl}
                   q_{\Delta} :  \J_{2}(\C) & \longrightarrow & \Z/2\Z \\
                    j  & \longmapsto       & Q(\Delta)+Q(j+\Delta) \end{array}  \right)  \end{array}$$
                    is an isomorphism of symplectic torsors.        \label{s}
 In particular it preserves the parity.
 \end{prop}

\section{Proof of  theorem \ref{f}}

          \subsection{Orbits of the action of $\mathrm{Sp}_{12}(\mathbb{Z}/2\mathbb{Z})$
          on the  sets of four distincts even theta characteristics}
           Let  $\C$ be a curve of genus $6$ over $\mathbb{C}$. For each set $\{\Delta_{i}\}_{i=1\dots4}$ of four distincts even
           theta characteristics on $\C$, we can associate the set $\{\delta_{i}\}_{i=1..4}$ of theta characteristics on $\C$, defined as
             $$  \delta_{i} \equiv 2K-\sum_{j \neq i} \Delta_{j} .     $$

         \begin{prop}
              Let  $C$ be a curve of genus $6$ over $\mathbb{C}$. The action of $\mathrm{Sp}_{12}(\mathbb{Z}/2\mathbb{Z})$
              on  the sets of four even theta characteristics has four orbits $A_{1},\dots,A_{4}$, given respectively by the following
              conditions :
                       \begin{enumerate}
                           \item  $\sum_{i=1}^{4} \Delta_{i} \equiv 2\K$.
                           \item   $\{\Delta_{i}\}_{i=1\dots4}$ is not in $A_{1}$ and all the $\delta_{i}$'s are even.
                           \item Exactly two of the  $\delta_{i}$'s are even.
                            \item All the $\delta_{i}$'s are odd.
                       \end{enumerate}
                    \label{orb}
         \end{prop}

              By Proposition \ref{s}, we can study the action of  $\mathrm{Sp}_{12}(\mathbb{Z}/2\mathbb{Z})$
          on the sets of four distinct even quadratic forms on  $\J_{2}(\C)$ :

                          \begin{lemm}  Let  $C$ be a curve of genus $6$ over $\mathbb{C}$.
                          Let $\{q, q+a_{1},q+a_{2},q+a_{3}\}$ and   $\{q, q+a_{1}^{'},q+a_{2}^{'},q+a_{3}^{'}\}$ be two
                          sets of four distinct even quadratic forms on $\J_{2}(\C)$ so that
                                                   $$\dim<a_{1},a_{2},a_{3}>
                                                                                =\dim<a_{1}^{'},a_{2}^{'},a_{3}^{'}>.$$
                            If  it exists $\sigma \in \mathfrak{S}_{3}$ so that
                                        $$\forall k,l \in \{1,2,3\} \ , \   a_{k}.a_{l}=a_{\sigma(k)}^{'}.a_{\sigma(l)}^{'},$$
                             then these two sets are in the same orbit of the action of $\mathrm{Sp}_{12}(\mathbb{Z}/2\mathbb{Z})$
                              on the sets of four distinct even quadratic forms on  $\J_{2}(\C)$.
                              \label{witt}
                 \end{lemm}

    \begin{description}
            \item[Proof : ]
                By taking into account  the parity conditions on the quadratic forms, we have :
                $$\begin{array}{c}
                       q(a_{1})=q(a_{1}^{'})      \\
                             q(a_{2})=q(a_{2}^{'}) \\
                             q(a_{3})=q(a_{3}^{'}).  \end{array}$$
               Moreover, $\dim<a_{1},a_{2},a_{3}>=\dim<a_{1}^{'},a_{2}^{'},a_{3}^{'}>$ and it exists $\sigma \in \mathfrak{S}_{3}$ so that
                                        $$\forall k,l \in \{1,2,3\}  , \   a_{k}.a_{k}=a_{\sigma(k)}^{'}.a_{\sigma(l)}^{'}.$$
               By Witt's theorem  (see for instance [\textbf{Sc}]) these conditions are equivalent to :
                $$ \exists  f \in \mathrm{O}(q) \  so \ that \ \forall k \in \{1,2,3\}  ,\  f(a_{k})=a_{\sigma(k)}^{'}. $$
                Finally, this is equivalent to the existence of $f$ in $\mathrm{O}(q)$ so that
                                            $$\begin{array}{l}
                           q_{1}=q_{i_{1}}^{'} \circ f^{-1}     \\
                            q_{2}=q_{i_{2}}^{'} \circ f^{-1}  \\
                             q_{3}=q_{i_{3}}^{'} \circ f^{-1}. \ \square  \end{array}   $$

\end{description}

 \begin{description}
   \item [Proof of Proposition \ref{orb} : ]
   As the action preserves the parity, each orbit is contained in one of the $A_{i}$'s.
   By this fact, we have to show the transitivity of this action
   on these sets : \\
Let  $\Delta_{1},\Delta_{2},\Delta_{3}$ and $\Delta_{4}$ be four even  theta characteristics on $C$. Let
$q_{1},q_{2},q_{3}$ and $q_{4}$ be the quadratic forms associated to these theta characteristics.
 For each $i$ in $\{1,2,3\}$, let $a_{i}$ be the element of  $\J_{2}(\C)$
 so that $$q_{i}=q_{4}+a_{i}.$$
   If $\ a_{1},a_{2}$ and $a_{3}$ are linearly dependent then
                                 $ \Delta_{4}+ \Delta_{1} \equiv \Delta_{2}+\Delta_{3}.$  \\
                                                If  $ \ a_{1},a_{2}$ and $a_{3}$ are
                          linearly independent then we have : \\ $\forall \{i,j\} \subseteq \{1,2,3\}$,
                               $$a_{i}.a_{j}=0 \Leftrightarrow  \Delta_{4}+\Delta_{i}-\Delta_{j} \ even.$$
 Thus, we have : \\
 $$\begin{array}{rll}
 \{\Delta_{i}\}_{i=1\dots4} \in A_{1} & \Leftrightarrow &  a_{1},a_{2} \ and \ a_{3} \ linearly \ dependent.\\
 \{\Delta_{i}\}_{i=1\dots4} \in A_{2}  &  \Leftrightarrow & \dim <a_{1},a_{2},a_{3}>=3 \ and \ \forall i \in \{1,2,3\} ,\ a_{i}.a_{j}=0 \\
 \{\Delta_{i}\}_{i=1\dots4}  \in A_{3}  &  \Leftrightarrow & \dim <a_{1},a_{2},a_{3}>=3 \ and  \\
   &  & \mathrm{Card} \left( \left\lbrace  \{i,j\} \subset \{1,2,3\} \mid a_{i}.a_{j}=1 \right\rbrace    \right) =1 \ or \ 2   \\
  \{\Delta_{i}\}_{i=1\dots4} \in A_{4} & \Leftrightarrow & \dim <a_{1},a_{2},a_{3}>=3 \ and \ \forall i \in \{1,2,3\} ,\ a_{i}.a_{j}=1 \end{array}$$

  Let  $q_{1}^{\ast},q_{2}^{\ast},q_{3}^{\ast}$ and $q_{4}^{\ast}$ be four even  quadratic forms on $\J_{2}(\C)$ corresponding to $\{\Delta_{i}^{\ast}\}_{i=1\dots4}$, an
  other set of  even theta characteristics on $C$. The transitivity of the action (on even quadratic forms) implies that there exists $M$ in $\mathrm{Sp}_{12}(\mathbb{Z}/2\mathbb{Z})$,
  so that $M.q_{4}^{\ast}=q_{4}$. If we put
  $$\begin{array}{c}
                           q_{1}^{'}:=M.q_{1}^{\ast}     \\
                            q_{2}^{'}:=M.q_{2}^{\ast}  \\
                             q_{3}^{'}:=M.q_{3}^{\ast},   \end{array}$$
                then $\{q_{1}^{'},q_{2}^{'},q_{3}^{'},q_{4}\}$ is in the same orbit as $\{q_{1}^{\ast},q_{2}^{\ast},q_{3}^{\ast},q_{4}^{\ast}\}$.

 By this fact, if  $\{\Delta_{i}\}_{i=1\dots4}$ and $\{\Delta_{i}^{\ast}\}_{i=1\dots4}$ are both in $A_{1},A_{2},A_{4}$, then the transitivity is
  obvious by lemma \ref{witt}. \\
  Now let us suppose that  $\{\Delta_{i}\}_{i=1\dots4}$ and $\{\Delta_{i}^{\ast}\}_{i=1\dots4}$ are in $A_{3}$.
  For each $i$ in $\{1,2,3\}$, let $a_{i}^{'}$ be the element of  $\J_{2}(\C)$
 so that $$q_{i}^{'}=q_{4}+a_{i}^{'}.$$
  If   $\mathrm{Card} \left( \left\lbrace  \{i,j\} \subset \{1,2,3\} \mid a_{i}.a_{j} \right\rbrace    \right)=
 \mathrm{Card} \left( \left\lbrace  \{i,j\} \subset \{1,2,3\} \mid a_{i}^{'}.a_{j}^{'} \right\rbrace    \right)$, then we can conclude
  again with lemma \ref{witt}. If not, we can suppose for instance, that
 $$\begin {array}{ccc}
      a_{1}.a_{2}=1 &  &  a_{1}^{'}.a_{2}^{'}=1  \\
      a_{1}.a_{3}=0 & and & a_{1}^{'}.a_{3}^{'}=1     \\
      a_{2}.a_{3}=0 &  & a_{2}^{'}.a_{3}^{'}=0.     \end{array} $$
By lemma \ref{witt},  there exist $a_{2}^{''}$ and $a_{3}^{''}$ in $\J_{2}(\C)$ so that
  $$\begin {array}{c}
      a_{2}.a_{2}^{''}=1   \\
      a_{2}.a_{3}^{''}=1      \\
      a_{2}^{''}.a_{3}^{''}=0,    \end{array} $$
    and $\{q_{4}+a_{2},q_{4}+a_{2}^{''},q_{4}+a_{3}^{''},q_{4}\}$
      is in the same orbit as
 $\{q_{4}+a_{1}^{'},q_{4}+a_{2}^{'},q_{4}+a_{3}^{'},q_{4}\}=\{q_{1}^{'},q_{2}^{'},q_{3}^{'},q_{4}\}$. If we put $$q:=q_{4}+a_{2},$$
then $$\begin{array}{l}
            \{q_{1},q_{2},q_{3},q_{4}\}=\{q+(a_{1}+a_{2}),q,q+(a_{3}+a_{2}),q+a_{2}\} \\
             \{q_{4}+a_{2},q_{4}+a_{2}^{''},q_{4}+a_{3}^{''},q_{4}\}=\{q,q+(a_{2}+a_{2}^{''}),q+(a_{2}+a_{3}^{''}),q+a_{2}\}. \end{array}$$
We can verify that
     $$\begin {array}{ccc}
      (a_{1}+a_{2}).(a_{2}+a_{3})=1 &  &  a_{2}.(a_{2}+a_{2}^{''})=1  \\
      (a_{1}+a_{2}).a_{2}=1 & and & a_{2}.(a_{2}+a_{3}^{''})=1     \\
      (a_{2}+a_{3}).a_{2}=0 &  & (a_{2}+a_{2}^{''}).(a_{2}+a_{3}^{''})=0.     \end{array} $$
 and conclude with the lemma \ref{witt}. $\square$
 \end{description}

\subsection{Theta characteristics  on  bi-elliptic curves}
  Let  $\C$ be a bi-elliptic curve of genus $6$. By definition this means that there exists an elliptic curve  $\mathrm{X}$ and
  a degree $2$ morphism :
    $$\pi:\C \longrightarrow \mathrm{X}$$
   The Riemann-Hurwitz theorem implies that  $\pi$ has $10$ ramification points denoted $R_{1},\dots,R_{10}$. \\
    \begin{lemm}([\textbf{A2}])
       A general bi-elliptic curve has fourty even, effective theta characteristics of the form
      $$R_{i}+\pi^{\ast}(D_{i})$$
       where  $R_{i}$ is one of the ramification points of $\pi$ and  $D_{i}$ is a degree two divisor on  $\mathrm{X}$.
  \end{lemm}

  \begin{rema}
    \begin{itemize}
     \item $R_{i}$ is the fixed point of the linear system $\mid R_{i}+\pi^{\ast}(D_{i} \mid$.
      \item $\J(\mathrm{X})$ has $3$ non zero points of order $2$ denoted by $F_{1},F_{2}$ and $F_{3}$.
     If for each $i$ in  $\{1,\dots,10\}$, we choose a divisor  $D_{i}$ so that  \\ $R_{i}+\pi^{\ast}D_{i}$ is
     an even theta characteristic on   $\C$, then the three other even, effective theta characteristics with fixed point   $R_{i}$ will be
                                 $$R_{i}+\pi^{\ast}(D_{i}+F_{j}) \ \ \ \ \ \ \ j=1\dots 3 $$
  \end{itemize}
  \end{rema}

  \begin{lemm}
       Let   $\Delta_{1}=R_{1}+\pi^{\ast}(D_{1})$, $\Delta_{2}=R_{2}+\pi^{\ast}(D_{2})$ and \\
        $\Delta_{3}=R_{3}+\pi^{\ast}(D_{3})$ be three even effective theta characteristics on a general
        bi-elliptic curve $\C$. One has
        $$\Delta_{1} +\Delta_{2} -\Delta_{3} \  even \ \Leftrightarrow \left\lbrace
         \begin{array}{l}  two \  at \ least \ of \\
                                     the\ R_{i} \ are \ the \ same     \end{array}  \right. $$  \label{bi}

  \end{lemm}

  \begin{description}
        \item[ Proof : ]
        First of all, let us notice that for each $j,k,l$ so that \\ $\{j,k,l\}=\{1,2,3\}$,
                      $$\Delta_{1}+\Delta_{2}-\Delta_{3} =R_{j}+R_{k}-R_{l}+\pi^{\ast}D_{j}
                      +\pi^{\ast}D_{k}-\pi^{\ast}D_{l} $$

            then $D_{j}+D_{k}-D_{l} \equiv E_{l}$ where $E_{l}$ is a degree two, effective divisor on $\mathrm{X}$ ;
            thus  $(\Leftarrow )$  is evident. \\
             Now, if we call  $P_{l}$, the point on  $\mathrm{X}$ so that   $\pi^{\ast}P_{l}=2R_{l}$, as
             on an elliptic curve every linear system of degree $d$ is of dimension  $d$, there exists $Q_{l}$,
             a point on   $\mathrm{X}$ so that
                      $P_{l}+Q_{l} \equiv E_{l}$ ;
            thus $$\begin{array}{cl} R_{j}+R_{k}-R_{l}+\pi^{\ast}D_{j}
                                                                                                   +\pi^{\ast}D_{k}-\pi^{\ast}D_{l}
                                                                                                   & \equiv   R_{j}+R_{k}-R_{l} +\pi^{\ast}{E_{l}} \\
                                                                               &\equiv   R_{j}+R_{k}-R_{l} +\pi^{\ast}{P_{l}}+\pi^{\ast}{Q_{l}} \\
                                                                                &\equiv   R_{1}+R_{2}+R_{3} +\pi^{\ast}{Q_{l}}
                                              \end{array}$$
              But then, by lemma \ref{bi}, for a general bi-elliptic curve, if we are considering, for instance, that $ R_{1}$ is the fixed point of the even
              theta characteristic
              $\Delta_{1}+\Delta_{2}-\Delta_{3}$, then $R_{2}+R_{3} +\pi^{\ast}{Q_{l}}$
               must be the pullback of an effective degree two divisor. This is possible if and only if  $R_{2}=R_{3}$. $\square$

  \end{description}

With this lemma, one checks easily the next result which ends our proof.

  \begin{prop}
      Let  $\C$ be a smooth  curve of genus $6$ with a degree two morphism
       $$\pi:\C \longrightarrow \mathrm{X}$$
       onto an elliptic curve $\mathrm{X}$.
       Let $R_{1},\dots,R_{10}$ be the ramification points of this morphism and let $D_{1},\dots,D_{10}$
       be some degree two effective divisors on  $\mathrm{X}$ so that
      $$\{R_{i}+\pi^{\ast}D_{i} \ \mid \ i \in \{1,\dots,10\}\}$$
        is a set of the even, effective theta characteristic on  $\C$. The following sets are in the four orbits of the action of
        $\mathrm{Sp}_{12}(\mathbb{Z}/2\mathbb{Z})$ on the sets of four distinct even theta characteristics on $\C$ :
                                    $$\begin{array}{l}
                                     \{R_{1}+\pi^{\ast}D_{1},R_{2}+\pi^{\ast}(D_{2}+F_{1}),R_{1}+\pi^{\ast}(D_{1}+
                       F_{1}) , R_{2}+\pi^{\ast}D_{2}\}  \in A_{1} \\
                            \{R_{1}+\pi^{\ast}D_{1},R_{2}+\pi^{\ast}(D_{2}+F_{2}),R_{1}+\pi^{\ast}(D_{1}+
                       F_{1}) , R_{2}+\pi^{\ast}D_{2} \}  \in A_{2} \\
                               \{R_{1}+\pi^{\ast}D_{1},R_{2}+\pi^{\ast}(D_{2}),R_{3}+\pi^{\ast}(D_{3}) , R_{1}+\pi^{\ast}(D_{2}+F_{1})\} \in A_{3}   \\
                    \{R_{1}+\pi^{\ast}D_{1},R_{2}+\pi^{\ast}D_{2},R_{3}+\pi^{\ast}D_{3} , R_{4}+\pi^{\ast}D_{4}\} \in A_{4}.
                  \end{array}$$

  \end{prop}

  \section{Proof of theorem \ref{trans} }
   First of all, we need to recall a result given by Teixidor I Bigas in  [\textbf{T}] : \\
   Let $\St_{1}$ be the scheme wich parametrizes the pairs $(\C,\Delta)$, where $\C$ is curve of
   genus $g$ on $\mathbb{C}$ and $\Delta$ is a theta characteristic of projective dimension $1$  on $\C$.
   For each $t=(\C,\Delta)$ in $\St_{1}$, there is an injective map :
    $$f:T_{\St_{1}(t)} \longrightarrow  H^{1}(T_{\C})$$
    As by Serre duality $ H^{1}(T_{\C})$ is dual to $H^{0}(2\K)$, we have :
       \begin{lemm}(Teixidor I Bigas)
              If  $F$ is the fixed part of the linear system $| \Delta |$ and  $R$ the ramification divisor of the corresponding morphism
         $\varphi_{\Delta}:\C \longrightarrow \PP^{1}$ then $\mathrm{Im}(f)=(\omega)^{\bot}$, where $\omega$
         is an element of  $H^{0}(2\K)$ with divisor  $R+2F$.
   \end{lemm}

    By using the particular expression of the theta characteristics on an
    \\ hyperelliptic curve (see [\textbf{A,C,G,H}] p 288),
     this lemma has the following consequence :
     Let $\C$ be an hyperelliptic curve of genus $g\geqslant3$,  $W:=\{p_{1},\dots,p_{2g+2}\}$ the set of its Weierstrass points, $H$ its hyperelliptic
           divisor, then for each even theta characteristic $\Delta$ on $\C$ of dimension $1$,  $$\Delta=H+p_{i_{1}}+\dots+p_{i_{g-3}}.$$
         $\mathrm{Im}(f)$ is the orthogonal of $\omega$ so that
                                $$\mathrm{div}(\omega)=p_{1}+\dots+p_{2g+2}+2(p_{i_{1}}+\dots+p_{i_{g-3}})$$
        \begin{prop}
         Let $(\C,\sigma)$ be an hyperelliptic element of $\mathcal{M}_{g}(2)$ ($g \geqslant 3$). Let $H$ be
         the hyperelliptic divisor on $C$. We choose $p_{1},\dots,p_{g-2}$, $g-2$ Weierstrass points on $\C$ and
         let $E$ be the divisor $p_{1}+\dots+p_{g-2}$. Then the $g-2$ thetanull divisors associated to the theta characteristics
         $\{H+E-p_{1},\dots,H+E-p_{g-2}\}$ intersect transversally at $(\C,\sigma)$.  \label{fin}
       \end{prop}

       \begin{description}
            \item[Proof : ]
            By the last remark, we have to prove that the linear subsystem of  $H^{0}(2\K)$ generated by
             $$\{ F_{k}=p_{1}+\dots+p_{2g+2}+2(E-p_{k})  \mid  k=1,\dots,g-2 \}$$
             has rank  $g-2$. Let us prove that  for each $k$ in $\{2,\dots,g-2\}$, $F_{k}$
             is not in the linear system generated by  $F_{1},\dots,F_{k-1}$.
             Let $\C$ be the Riemann surface
                             $$y^2=\prod_{i=1}^{2g+2}(x-x_{i}) \ \ ,\ \ x:\C \longrightarrow \PP^1$$
            so that for each  $i$ in  $\{1,\dots,2g+2\}$ , $x(p_{i})=x_{i}$ .
               We have then to verify that there does not exist  $(\la_{1}:\dots:\la_{k-1})$ in $\PP^{k-2} $ so that
                $2(E-p_{k})$ is the divisor  $f^{\star} 0$ where  $$f=\la_{1}( x-x_{2})\dots
                             (x-x_{g-2})+\dots+\la_{k-1}(x-x_{1})\dots(x-x_{k-2})(x-x_{k})\dots(x-x_{g-2})$$
                Let us suppose it is false. As  $2(E-p_{k})$ is  $2p_{1}+\dots+2p_{k-1}+ 2p_{k+1}+  \dots+ 2p_{g-2} $ this implies
                             $$f(x_{1})=\la_{1}(x_{1}-x_{2})( x_{1}-x_{2})\dots(x_{1}-x_{g-2})=0 \ \Rightarrow \ \la_{1}=0$$
                             then
                             $$f(x_{2})=0 \ \Rightarrow \ \la_{2}=0$$
                             then
                              $$\vdots$$
                              $$f(x_{k-1})=\la_{k-1}(x_{k-1}-x_{1})\dots( x_{k-1}-x_{k-2})( x_{k-1}-x_{k})\dots(x_{k-1}-x_{g-2})=0
                             $$ $$\ \Rightarrow  \la_{k-1}=0$$
                             Finally  we would have $\la_{1}=\dots=\la_{k-1}=0$ which it is absurd. $\square$

          \end{description}
Each connected component $Y$ of  $\mathcal{H}_{g}(2)$ is uniquely determined by a symplectic isomorphism
$$c_{Y}: (\Z/2\Z)^{2g} \longrightarrow  P_{2}^{+}(\{1,\dots,2g+2\})$$
                    and an isomorphism between symplectic torsors
                      $$\widetilde{c_{Y}}: (\Z/2\Z)^{2g} \longrightarrow  T(\{1,\dots,2g+2\})$$
 in the following way : \\
  Let $U$ be the open subvariety of $\mathbb{C}^{2g+2}$ consisting of points with distincts coordinates. To each point
        $\xi=(\xi_{1},\dots,\xi_{2g+2})$ in $U$, we associate the hyperelliptic curve $\C(\xi)$ with Weierstrass
        points $\xi_{1},\dots,\xi_{2g+2}$. We denote by $W(\xi)$ this set.
               For each $\xi$ in $U$, the bijection $[1,\dots,2g+2] \longrightarrow [\xi_{1},\dots,\xi_{2g+2}]$ induces
          a symplectic isomorphism
            $$ \alpha_{\xi}:  P_{2}^{+}(\{1,\dots,2g+2\})   \longrightarrow   P_{2}^{+}(W(\xi)) $$
          and an isomorphism between symplectic torsors :
                 $$ \widetilde{\alpha_{\xi}}:  T(\{1,\dots,2g+2\})   \longrightarrow   T(W(\xi)) .$$
  $Y$ is the subspace of $\mathcal{M}_{g}(2)$ defined as
              $$\{ (\C(\xi),\sigma_{\xi}) \in  \mathcal{H}_{g,2} \mid \xi \in U\},$$
 so that for each $\xi$ in $U$, the following diagrams
   $$ \xymatrix{    P_{2}^{+}(\{1,\dots,2g+2\}) \ar[r]^-{\alpha_{\xi}} &  P_{2}^{+}(C(\xi)) \ar[rr]^-{\beta_{\C(\xi)}}
                    &       &      \J_{2}(\C(\xi))   \\  &  \ar[ul]_{c_{Y}} \ar[urr]_{\sigma_{\xi}} (\Z/2\Z)^{2g} & & } $$

                           $$ \xymatrix{    T(\{1,\dots,2g+2\}) \ar[r]^-{\widetilde{\alpha_{\xi}}} &  T(W(\xi)) \ar[rr]^-{\widetilde{\beta_{\C(\xi)}}}
                    &       &       \mathcal{S}(\C(\xi))   \\  & \ar[urr]_{\widetilde{\sigma_{\xi}}} \ar[ul]_{\widetilde{c_{Y}}} (\Z/2\Z)^{2g} & & } $$
 are commutative.
Now let $\xi_{0}$ be in $U$. Let us choose $g-2$ thetanull divisors  which correspond on $\C(\xi_{0})$, to $g-2$
theta characteristics in the configuration of Proposition \ref{fin}. By the last diagram, one sees that for each
 $\xi$ in $U$, these  $g-2$ thetanull divisors correspond to $g-2$ theta characteristics
              on $\C(\xi)$  in the same configuration as on $\C(\xi_{0})$. By this fact, $Y$
           is  a component of the intersection of these $g-2$ thetanull divisors. $\square$

 \section*{Acknowledgment}  I would like to thank my thesis advisor A. Beauville for making useful suggestions that lead to the
 results of this paper.

\end{document}